\documentclass[11pt]{article}

\usepackage{inputenc}
\usepackage{makeidx}
\usepackage{amssymb}
\usepackage{amsfonts}
\usepackage{amsmath}
\usepackage{euscript}
\usepackage{theorem}
\usepackage{enumerate}
\usepackage{graphics}

\theorembodyfont{\rm}      

\theoremstyle{change}      

\begingroup

\newtheorem{thm}{Theorem\hskip 5mm}[section]

\newtheorem{cor}[thm]{Corollary\hskip 5mm}

\newtheorem{lem}[thm]{Lemma\hskip 5mm}

\newtheorem{note}[thm]{Note\hskip 5mm}

\endgroup

\begingroup

\endgroup

\begingroup

\endgroup

\def\Z{{\mathbb Z}}
\def\C{{\mathbb C}}
\def\R{{\mathbb R}}
\def\m{{\mathfrak m}}
\def\ker{\mathrm{ ker\,}}

\begin{document}

\begin{center}
{\bf \large SUBALGEBRAS OF MATRIX ALGEBRAS\\ GENERATED BY
COMPANION MATRICES}
\end{center}

\begin{center}{{\large {\sc N. H. Guersenzvaig}}\\
{\small A\!v. Corrientes 3985 6A, (1194) Buenos Aires, Argentina\\ email: nguersenz@fibertel.com.ar}}
\end{center}

\begin{center}
and
\end{center}

\begin{center}{{\large {\sc Fernando Szechtman}}\\
{\small Department of Mathematics and Statistics, University of Regina, Saskatchewan, Canada\\
email: fernando.szechtman@gmail.com}}
\end{center}


\begin{abstract}
Let $f,g\in \Z[X]$ be monic polynomials of degree $n$ and let 
$C,D\in M_n(\Z)$ be the corresponding companion matrices. We find necessary
and sufficient conditions for the subalgebra $\Z\langle C,D\rangle$ to be a sublattice
of finite index in the full integral lattice $M_n(\Z)$, in which case we compute the exact value of this index 
in terms of the resultant of $f$ and $g$. If $R$ is a commutative ring
with identity we determine when $R\langle C,D\rangle=M_n(R)$, in which case a presentation for 
$M_n(R)$ in terms of $C$ and $D$ is given.

\end{abstract}



\smallskip


\section{Introduction} About twenty years ago a question of Chatters [C1]
generated a series of articles concerned with the problem of
identifying full matrix rings. We refer the reader to the papers
[A], [AMR], [C2], [LRS], [R] cited in the bibliography for more
details. In particular, very simple presentations of full matrix
rings, involving just two generators, were obtained.

In this paper we concentrate on the algebra generated two matrices
$A,B\in M_n(R)$, where $R$ is a commutative ring with identity and
$n\geq 2$. Is it possible to find a presentation for $R\langle
A,B\rangle$? If $A$ and $B$ happen not to generate $M_n(R)$, can
we somehow measure the degree of this failure? Adopting a more
precise and geometric viewpoint, we look at $M_n(\Z)$ as an
integral lattice in $M_n(\R)$ and ask when will the sublattice
$\Z\langle A,B\rangle$ have maximal rank and, in that case, what
will be its index in the full lattice $M_n(\Z)$. The answers to
these questions depend on more specific information about $A$ and
$B$.

Focusing attention on two companion matrices $C,D\in M_n(R)$ of
monic polynomials $f,g\in R[X]$ of degree $n$, section \ref{gen}
gives necessary and sufficient conditions for $C$ and $D$ to
generate $M_n(R)$, while section \ref{presmn} determines how they
do it. If $R$ is a unique factorization domain, section
\ref{prescd} exhibits a presentation of $R\langle C,D\rangle$,
proves it to be a free $R$-module, and computes its rank.

In section \ref{indic} we find the exact index of $\Z\langle
C,D\rangle$ in $M_n(\Z)$ and extend this result to other number
rings. The index is obtained by means of a determinantal identity,
found in section \ref{D}, which is of independent interest and
valid under no restrictions on $R$.

We will keep the above notation as well as the following. Let
$R[X,Y]$ be the $R$-span of $X^i Y^j$ in $R\langle X,Y\rangle$,
where $0\leq i,j$. We have a natural map $R\langle X,Y\rangle\to
M_n(R)$ sending $X$ to $A$ and $Y$ to $B$. Let $R[A,B]$ stand for
the image of $R[X,Y]$ under this map. Since $A$ and $B$ are
annihilated by their characteristic polynomials, we see that
$R[A,B]$ is spanned by $A^iB^j$, where $0\leq i,j\leq n-1$.
Clearly $R[A,B]\subseteq R\langle A,B\rangle$, with equality if
and only if $R[A,B]$ is a subalgebra, which is definitely not
always true. Perhaps surprisingly, section~\ref{cesd} proves that
$R[C,D]=R\langle C, D\rangle$. A more detailed discussion of this
is given in section \ref{pj}.

The resultant of $f$ and $g$ will be denoted by $R(f,g)$. A fact
used repeatedly below is that $R(f,g)$ is a unit if and only if
$f$ and $g$ are relatively prime when reduced modulo every maximal
ideal of~$R$.

\section{A theorem of Burnside}

For the record, we state here general conditions for a subset $S$
of $M_n(R)$ to generate $M_n(R)$ as an algebra. The field case
follows from Burnside's Theorem (see $\S 27$ of [CR]) whereby one obtains the general
case by localization.

\begin{thm}
\label{dense} Let $F$ be a field and
 let $S$ be subset of $M_n(F)$. Then the subalgebra
generated by $S$ is the full matrix algebra $M_n(F)$ if and only
if the following conditions~hold:

\begin{enumerate}
\item [(C1)] The only matrices in $M_n(F)$ commuting with all
matrices in $S$ are the scalar matrices.

\item [(C2) ] The only subspaces of the column space $V=F^n$ that
are invariant under the action of all matrices in $S$ are 0 and
$V$.
\end{enumerate}
\end{thm}

\begin{note} In the above notation, if $F$ is a subfield of $K$ then
$\mathrm{dim}_F F\langle S\rangle=\mathrm{dim}_K K\langle S\rangle$,
so  $F\langle S\rangle=M_n(F)$ if and only if $K\langle S\rangle=M_n(K)$.
\end{note}

\begin{thm}
\label{local} For each maximal ideal $\m$ of $R$, let $\Lambda_\m:
M_n(R)\to M_n(R/\m)$ be the ring epimorphism associated to the
projection $R\to R/\m$.  Let $S$ be subset of $M_n(R)$. Then
$$
R\langle S\rangle=M_n(R)\Leftrightarrow (R/\m)\langle
\Lambda_\m(S)\rangle=M_n(R/\m)
$$
for every maximal ideal $\m$ of $R$.
\end{thm}

\noindent{\sl Proof.} One implication is obvious. For the other,
suppose that $(R/\m) \langle\Lambda_\m(S)\rangle=M_n(R/\m)$ for
all maximal ideals $\m$ of $R$. This is equivalent to
$\Lambda_\m(R\langle S\rangle)=\Lambda_\m(M_n(R))$, that is,
$R\langle S\rangle+\ker \Lambda_\m=M_n(R)$, which obviously means,
$R\langle S\rangle+\m\, M_n(R)=M_n(R)$, or just $\m
\big(M_n(R)/R\langle S\rangle\big) = M_n(R)/R\langle S\rangle$, for
all maximal ideals $\m$ of $R$. The quotient, say
$U=M_n(R)/R\langle S\rangle$, in this last statement is an
$R$-module.

Now $M_n(R)$ is a finitely generated $R$-module, and hence so is
$U$. We are thus faced with a finitely generated $R$-module,
namely $U$, such that $\m U=U$ for all maximal ideals $\m$ of $R$.
Localizing $R$ and $U$ at $\m$ (see chapter 3 of [AM]), we obtain
that ${\mathfrak M} U_\m=U_\m$, where ${\mathfrak M}$ is the
maximal ideal of $R_\m$. As $U_\m$ is a finitely generated
$R_\m$-module, it now follows from Nakayama's Lemma that $U_\m=0$
for all maximal ideals $\m$ of $R$. As all localizations at
maximal ideals of $U$ are zero, it follows from Proposition 3.8 of
[AM] that $U$ itself is zero, that is, $R\langle S\rangle=M_n(R)$. \quad $\blacksquare$

\section{Matrices commuting with $C$ and $D$}

We fix the following notation for the remainder of the paper:
$e_1, \dots , e_n$ will stand for the canonical basis of the
column space $R^n$ and $R_n[X]$ for the $R$-submodule of $R[X]$
with basis $1,X,\dots,X^{n-1}$. If $p\in R_n[X]$ then $[p]$ stands
for the coordinates of $p$ relative to this basis. Recall that $C$
is the companion matrix to $f=f_0+f_1X+\cdots+f_{n-1}X^{n-1}+X^n$,
that is
$$
C=\left(%
\begin{array}{ccccc}
  0 & 0 & \cdots & 0 & -f_0 \\
  1 & 0 & \cdots & 0 & -f_1 \\
  0 & 1 & \cdots & 0 & -f_2 \\
  \vdots & \vdots & \cdots & \vdots & \vdots \\
  0 & 0 & \cdots & 1 & -f_{n-1} \\
\end{array}%
\right).
$$
It is an easy exercise to verify that, as in the field case, the
minimal polynomial of $C$ is $f$. Thus $I,C, \dots , C^{n-1}$ is
an $R$-basis of $R[C]$. If $A\in R[C]$ we write $[A]$ for the
coordinates of $A$ relative to this basis. The next result is borrowed from [GS].

\begin{lem}
\label{coord}
  If $A\in R[C]$ then $A = ([A] \;C[A] \;\dots
\;C^{n-1}[A])$.
\end{lem}

\noindent{\sl Proof.} We have $A = y_0I + y_1C + \cdots +
y_{n-1}C^{n-1}$  with $y_j\in R$. Multiplying both sides by $e_1$
gives $Ae_1 = y_0e_1 + y_1e_2 + \cdots + y_{n-1}e_n = [A]$.  If
$2\le j\le n$ then $Ae_j = AC^{j-1}e_1= C^{j-1}Ae_1 = C^{j-1}[A]$.
Thus the matrices in question have the same columns. $\quad
\blacksquare$

\begin{lem}
\label{columnspq}
 Let $p,q\in R_n[X]$. Then $p(C)[q]=q(C)[p]$. Also,
$p(C)[q]=0\Leftrightarrow f|pq$.
\end{lem}

\noindent{\sl Proof.} The first column of $p(C)q(C)=q(C)p(C)$
equals both $p(C)q(C)e_1=p(C)[q]$ and $q(C)p(C)e_1=q(C)[p]$ by
Lemma \ref{coord}. The remaining columns of $pq(C)$ are
$p(C)q(C)e_i=p(C)C^i[q]=C^ip(C)[q]$, $1\leq i\leq n-1$, so
$p(C)[q]=0 \Leftrightarrow pq(C)=0 \Leftrightarrow f|pq$. $\quad
\blacksquare$

\begin{lem}
\label{esescalar} Let $A\in R[C]$. Assume $A_{n,1}=\cdots=
A_{n,n-1}=0$. Then $A$ is scalar.
\end{lem}

\smallskip
\noindent {\sl Proof.} By Lemma \ref{coord} and hypothesis we have
\begin{equation}
\label{holaR} A=A_{1,1}I+A_{2,1}C+\cdots+A_{n-1,1}C^{n-2}.
\end{equation}
If $n=2$ we are done. Otherwise, applying both sides to $e_2$
gives $$Ae_2=A_{1,1}e_2+A_{2,1}e_3+\cdots+A_{n-1,1}e_n.$$ By
hypothesis $e_n$ does not appear in the second column of $A$,
namely $Ae_2$. Therefore $A_{n-1,1}=0$. Going back to
(\ref{holaR}), eliminating $A_{n-1,1}C^{n-2}$, and repeating the
argument with $e_3,\dots,e_{n-1}$ yields $A=A_{1,1}I$, as
required. \quad $\blacksquare$

\begin{lem}
\label{inter}
 Suppose that $R$ is an integral domain and that
$f\neq g$. Then the only matrices in $M_n(R)$ that commute with
$C$ and $D$ are the scalar matrices.
\end{lem}

\noindent{\sl Proof.} Suppose $A\in M_n(R)$ commutes with $C$ and
$D$. Then $A$ commutes with $Z=C-D$. But the first $n-1$ columns
of $Z$ are equal to zero and, by hypothesis, at least one entry of
the last column of $Z$ is not zero. Applying these facts to the
equation $AZ=ZA$ immediately gives $A_{n,1}=\cdots= A_{n,n-1}=0$.
Thus $A$ is scalar by Lemma \ref{esescalar}. \quad $\blacksquare$

\section{Common invariant subspaces under companion matrices}

The invariant subspaces under a single cyclic transformation are well-known and easily
determined. We gather all relevant information below.

\begin{lem}
\label{subespacio}
 Let $F$ be a field and $V$ a vector space over  $F$ of finite dimension $n$.
Let $T:V\to V$ be a cyclic linear transformation with cyclic
vector $v$ and minimal polynomial~$f$. Then the distinct
$T$-invariant subspaces of $V$ are of the form
$$
V(g)=V(g,T)=\{g(T)x\,|\, x\in V\},
$$
where $g$ runs through the monic factors of $f$. Moreover, $V(g)$
has dimension $n-\mathrm{deg}\,g$ and the $T$-conductor of $V$
into $V(g)$ is precisely $g$.
\end{lem}

\begin{lem}
\label{coprimo} Let $F$ be a field.
 Let $f_1, \dots , f_m$ be monic polynomials in $F[X]$ of
degree $n$. Then their companion matrices $C_{f_1}, \dots ,
C_{f_m}$ have a common invariant subspace different from 0 and
$V=F^n$ if and only $f_1, \dots , f_m$ have a common monic factor
whose degree is strictly between 0 and $n$.
\end{lem}

\noindent{\sl Proof.} Suppose $h$ is a common monic factor to all
$f_1, \dots ,f_m$ of degree strictly between~0 and $n$. By Lemma
\ref{subespacio} if $1\leq i\leq m$ then $V(h,C_{f_i})$ is a
$C_{f_i}$-invariant subspace of $V$ of dimension
$m=n-\mathrm{deg}\,h$, which is strictly between 0 and $n$. Now a
basis for $V(h,C_{f_i})$ is
$$
h(C_{f_i})e_1, C_{f_i}h(C_{f_i})e_1, \dots ,
C_{f_i}^{m-1}h(C_{f_i})e_1,
$$
which by Lemma \ref{coord} equals
$$
[h(X)], [Xh(X)], \dots , [X^{m-1}h(X)].
$$
Therefore all these subspaces are equal to each other.

Suppose conversely that $W$ is subspace of $V$ different from 0
and $V$ and invariant under $C_{f_1}, \dots , C_{f_m}$.  By Lemma
\ref{subespacio} we have $W=V(h_i,C_{f_i})$, where $h_i$ is a
monic factor of $f_i$ for each $i$. All $h_i$ have the same degree
and this degree is strictly between 0 and $n$, also by Lemma
\ref{subespacio}. We claim that the $h_i$ are all equal to $h_1$.
Indeed if $i>1$ then
$$[h_i]=h_i(C_{f_i})e_1\in V(h_i,C_{f_i})=W=V(h_1,C_{f_1}),$$
and therefore
$$[h_i]=s(C_{f_1})h_1(C_{f_1})e_1$$
for some $s\in F[X]$ of degree less than $n-\mathrm{deg}\, h_1$.
Hence by Lemma \ref{coord} $[h_i]=[sh_1]$ and therefore
$h_i=sh_1$. But $h_i$ and $h_1$ are monic of the same degree, so
$s=1$. \quad $\blacksquare$

\section{Generation of $M_n(R)$ by companion matrices}
\label{gen}

\begin{thm} Let $f_1, \dots , f_m$, $m\geq 2$, be monic polynomials in $R[X]$ of
degree $n$ with companion matrices $C_{f_1}, \dots , C_{f_m}$.
Then $R\langle C_{f_1},\dots, C_{f_m}\rangle=M_n(R)$ if and only
if $f_1,\dots,f_m$ are relatively prime when reduced modulo every
maximal ideal of $R$.
\end{thm}

\noindent{\sl Proof.} By Theorem \ref{local} $R\langle
C_{f_1},\dots, C_{f_m}\rangle=M_n(R)$ if and only if this equality
is preserved when $f_1,\dots,f_m$ and $R$ are reduced modulo every
maximal ideal. But at the field level, generation is equivalent to
the given polynomials being relatively prime, by Theorem
\ref{dense} and Lemmas \ref{inter} and \ref{coprimo}. \quad
$\blacksquare$

\begin{cor}
\label{unidad} $R\langle C,D\rangle=M_n(R)$ if and only if
$R(f,g)$ is a unit.
\end{cor}

\noindent{\bf Remark.} This does not generalize to arbitrary
matrices. Indeed, if $F$ is field then while two distinct Jordan
blocks in $M_n(F)$ have relatively prime minimal polynomials,
they share a common eigenvector, so they cannot generate the full
matrix algebra.

\section{The identity $R[C,D]=R\langle C,D\rangle=R[D,C]$}
\label{cesd}

\begin{lem}
\label{tecnico} Let $R\langle A,B\rangle$ be an $R$-algebra, where
$B$ is integral over $R$ of degree at most $n$. Then the following
three statements are equivalent:
\begin{enumerate}
\item[(a)] $B^jA\in R[A,B]$ for all $1\leq j\leq n-1$.

\item[(b)] $R[A,B]=R\langle A,B\rangle$.

\item[(c)] $(A-B)B^j(A-B)\in R[A,B]$ for all $0\leq j\leq n-2$.
\end{enumerate}
\end{lem}

\noindent{\sl Proof.} As $B$ is integral over $R$ of degree at
most $n$, condition (a) ensures that $R[A,B]$ is invariant under
right multiplication by $A$, which easily implies (b). On the
other hand, it is clear that (b) implies (c). Suppose finally that
(c) holds. We wish to prove that $B,BA,B^2A,\dots,B^{n-1}A$ are in
$R[A,B]$. We show this by induction. Clearly $B\in R[A,B]$.
Suppose $0<j\leq n-1$ and $B^{j-1}A\in R[A,B]$. By (c)
$$
B^{j+1}-B^jA-AB^j+AB^{j-1}A=(A-B)B^{j-1}(A-B)\in R[A,B].
$$
By definition $B^{j+1},AB^j\in R[A,B]$, while $A(B^{j-1}A)\in
R[A,B]$ by inductive hypothesis. Hence $B^jA\in R[A,B].\quad
\blacksquare$

\begin{lem}
\label{Z} Suppose the first $n-1$ columns of $Z\in M_n(R)$ are
equal to 0 and its last column has entries $z_1,\dots,z_n$. Let
$Q\in M_n(R)$ have entries $q_1,\dots,q_n$ in its last row. Then
$$
ZQZ=(q_1z_1+\cdots+q_nz_n)Z.
$$
\end{lem}

\smallskip
\noindent{\sl Proof.} We have
$$
\begin{aligned}
ZQZ &=\left(%
\begin{array}{cccc}
  0 & \cdots & 0 & z_1 \\
  \vdots &  & \vdots & \vdots \\
  0 & \cdots & 0 & z_n \\
\end{array}%
\right)\left(%
\begin{array}{ccc}
  * & \cdots  & * \\
  \vdots &  & \vdots \\
  q_1 & \cdots  & q_n \\
\end{array}%
\right)\left(%
\begin{array}{cccc}
  0 & \cdots & 0 & z_1 \\
  \vdots &  & \vdots & \vdots \\
  0 & \cdots & 0 & z_n \\
\end{array}%
\right)\\
&=\left(%
\begin{array}{ccc}
  z_1q_1 & \cdots & z_1q_n \\
  \vdots &  &  \vdots \\
   z_n q_1& \cdots  & z_n q_n\\
\end{array}%
\right)\left(%
\begin{array}{cccc}
  0 & \cdots & 0 & z_1 \\
  \vdots &  & \vdots & \vdots \\
  0 & \cdots & 0 & z_n \\
\end{array}%
\right)=(q_1z_1+\cdots+q_nz_n)Z.\quad \blacksquare
\end{aligned}
$$

\begin{cor}
\label{iguales} Suppose that $A,B\in M_n(R)$ share the first $n-1$
columns. Then $R[A,B]=R\langle A,B\rangle=R[B,A]$. In particular,
this holds when $A=C$ and $B=D$.
\end{cor}

\noindent{\sl Proof.} This follows at once from Lemmas
\ref{tecnico} and \ref{Z}.$\quad \blacksquare$

\bigskip \noindent{\bf Remark.} In general it is false that
$R\langle A, B\rangle=R[A,B]$ for arbitrary matrices $A$ and $B$,
even when $M_n(R)=R\langle A, B\rangle$. Indeed, consider the case
when $R=F$ is a field, $n\geq 3$, $A$ is a diagonal matrix with
distinct diagonal entries and $B$ is the all-ones matrix. The only
matrices commuting with $A$ must be diagonal and the only diagonal
matrices commuting with $B$ are scalar. Moreover, the only
non-zero subspaces of $V=F^n$ invariant under $A$ are spanned by
non-empty subsets of $e_1,\dots,e_n$ and none of them is
$B$-invariant except for $V$ itself. It follows from Burnside's
Theorem that $M_n(F)=F\langle A,B\rangle$. If we had $F\langle
A,B\rangle=F[A,B]$ then the $n^2$ matrices $A^iB^j$, with $0\leq
i,j\leq n-1$, spanning $F[A,B]$, would necessarily be linearly
independent, but they are not since $B^2=nB$.

\section{The polynomials $p_0,p_1,\dots,p_{n-1}\text{ behind }R[C,D]=R\langle C,D\rangle$}
\label{pj}

Since $R\langle C, D\rangle=R[C,D]$ or, equivalently, $R[C,D]$ is
invariant under right multiplication by $C$, there must exist
$n-1$ polynomials $P_1, \dots , P_{n-1}\in R[X,Y]$ satisfying:
\begin{equation}
\label{polyecu} D^jC=P_j(C,D),\quad j=1,\dots , n-1.
\end{equation}
In this section we define and explore an explicit sequence of
polynomials satisfying (\ref{polyecu}).

For the remainder of the paper we let $s=g-f\in R_n[X]$. Write
$a_{j}$ for the $(n,j)$-entry of $s(D)$ and set $u_n=e_n^t$. Using
first Lemma \ref{Z} and then Lemma \ref{coord} we see that
\begin{equation}
\label{aji} (C-D)D^{j-1}(C-D)\! = \! u_n D^{j-1}[s](C-D)=
u_ns(D)e_j(C-D)=a_j(C-D),\, 1\leq j\leq n.
\end{equation}

\begin{thm}
\label{polinomios} Define $p_0,p_1,\dots,p_{n-1}\in R_n[X]$ and
$P_0,P_1,\dots,P_{n-1}\in R[X,Y]$ by
$$p_0(X)=1,\quad
p_j(X) = X^j - a_{1}X^{j-1} - \cdots - a_{j-1}X - a_{j}, \quad
j=1, \dots , n-1,
$$
$$
P_j(X,Y) = p_j(X)(X-Y) + Y^{j+1}, \quad j=0, \dots , n-1.
$$
Then
\begin{enumerate}

\item [(a)]  $p_j(C)(C-D)=D^j(C-D)$ for all $0\leq j\leq n-1$.

\item[(b)] The polynomials $P_1,\dots,P_{n-1}\in R[X,Y]$ satisfy
(\ref{polyecu}).

\item[(c)] If $P=([p_0]\,[p_1]\dots [p_{n-1}])\in M_n(R)$ then
$g(C)P=-f(D)$.

\item[(d)] If $q_0, q_1, \dots , q_{n-1}\in R_n[X]$ and
$Q=([q_0]\,[q_1]\,\dots \,[q_{n-1}])\in M_n(R)$ then
$$D^j(C-D)\!=\!q_j(C)(C-D)\text{ for all }0\leq j\leq  n-1\Leftrightarrow
g(C)Q\!=\!-f(D)\Leftrightarrow g(C)(Q-P)\!=\!0.$$

\item[(e)] If $R(f,g)$ is a unit then $p_0,p_1,\dots,p_{n-1}$ is
the only sequence in $R_n[X]$ satisfying (a).

 \item [(f)]  If $s$ is a constant then $P_j(X,Y) = X^{j+1} + Y^{j+1}- X^jY$ for all $0\leq j\leq
n-1$.

\end{enumerate}
\end{thm}

\noindent{\sl Proof.} It is clear that $p_0(C)(C-D)=D^0(C-D)$. Let
$0< j\le n-1$ and suppose that $p_{j-1}(C)(C-D)=D^{j-1}(C-D)$.
Then (\ref{aji}) and the identity $p_j(X) = Xp_{j-1}(X) - a_{j}$
yield
\begin{align*}
p_j(C)(C-D) & = (Cp_{j-1}(C)-a_jI)(C-D)=
Cp_{j-1}(C)(C-D)-a_j(C-D)\\
&=CD^{j-1}(C-D)-(C-D)D^{j-1}(C-D)=D^j(C-D).
\end{align*}
This proves (a), which clearly implies (b). Note that
$q_j(C)(C-D)=D^j(C-D)$ can be written as $q_j(C)[s]=D^j[s]$, where
$0 \leq j\leq n-1$, which by Lemma \ref{columnspq} translates into
$s(C)Q=s(D)$, that is, $g(C)Q=-f(C)$. The sequence
$p_0,p_1,\dots,p_{n-1}$ does satisfy (a), so $(c)$ is true, whence
$g(C)Q=-f(C)\Leftrightarrow g(C)(Q-P)=0$, completing the proof of
(d). If $R(f,g)$ is a unit then $g(C)$ is invertible, in which
case $g(C)(Q-P)=0$ implies $Q=P$. This gives (e). If $s$ is a
constant then $a_j=0$ for all $1\leq j\leq n-1$, so $p_j=X^j$ and
a fortiori $P_j(X,Y)=X^j(X-Y)+Y^{j+1}=X^{j+1} + Y^{j+1}- X^jY$ for
all $0\leq j\leq n-1$.$\quad \blacksquare$

\section{A presentation of $M_n(R)$}
\label{presmn}

\begin{thm}
\label{thA} Suppose $R(f, g)$ is a unit.  Let $P_1,\dots, P_{n-1}$
be polynomials in $R[X,Y]$ as defined in Theorem \ref{polinomios}
or, more generally, be arbitrary as long as they satisfy
(\ref{polyecu}). Then the matrix algebra $M_n(R)$ has
presentation:
$$
\langle X,Y\,|\, f(X)=0,\, g(Y)=0,\, Y^jX=P_j(X,Y),\quad
j=1,\dots,n-1\rangle.
$$
In the particular case when $g-f$ is a unit in $R$, the matrix
algebra $M_n(R)$ has presentation
$$ \langle X,Y\,|\, f(X)=0,\, g(Y)=0,\,
Y^jX+X^jY=X^{j+1}+Y^{j+1},\quad j=1,\dots,n-1\rangle.
$$
\end{thm}

\noindent{\sl Proof.} Write $\Omega:R\langle X,Y\rangle\to
R\langle C,D\rangle$ for the natural $R$-algebra epimorphism that
sends $X$ to $C$ and $Y$ to $D$. Let $K$ be the kernel of
$\Omega$. Set $S=R\langle X,Y\rangle/K$ and let $A$ and $B$ be the
images of $X$ and $Y$ in $S$. We have $S=R[A,B]$ by Lemma
\ref{tecnico}, and it is clear that $S$ is $R$-spanned by
$A^iB^j$, $0\leq i,j\leq n-1$. If $t\in\ker \Omega$ then $t$ is a
linear combination of the $A^iB^j$. The images of these under
$\Omega$ are linearly independent, as $M_n(R)$ is free of rank
$n^2$ and, by Corollary \ref{unidad}, the $n^2$ matrices $C^iD^j$
span $M_n(R)$. Hence $t=0$.

If $g-f$ is a unit in $R$ then so is $R(f,g)$. Therefore the last
statement of the theorem follows from above and part (f) of
Theorem \ref{polinomios}. $\quad \blacksquare$

\medskip

As an illustration, let $R=\mathbb Q$, $f=X^n-2$, $g=X^n-3$. Let
$\alpha,\beta$ stand for the real $n$-th roots of 2 and 3,
respectively. Theorem \ref{thA} says that $M_n(\mathbb Q)=\mathbb
Q(\alpha)\mathbb Q(\beta)$, where $\mathbb Q(\alpha)$ and $\mathbb
Q(\beta)$ are embedded as maximal subfields of $M_n(\mathbb Q)$
which intersect only at $\mathbb Q$ and multiply according to the
rules:
$$\beta^j\alpha+\alpha^j\beta=\alpha^{j+1}+\beta^{j+1},\quad 1\leq j\leq
n-1.$$

\section{A presentation of $R\langle C,D\rangle$}
\label{prescd}

\begin{lem}
\label{h} Suppose that $d\in R[X]$ is a common monic factor of $f$
and $g$. Let $f=hd$, where $h\in R_n[X]$. Then $h(C)C=h(C)D$.
\end{lem}

\medskip
\noindent{\sl Proof.} By hypothesis $d|s$, whence $hd|hs$. But
$hd=f$, so Lemma \ref{columnspq} gives $h(C)[s]=0$, which implies
$h(C)(C-D)=(0,\dots,0,h(C)[s])=0$. $\quad \blacksquare$

\begin{thm}
\label{basisfg} Let $R$ be a unique factorization domain and let
$m=\mathrm{deg}\,\gcd(f,g)$. Then $R\langle C,D\rangle=R[C,D]$ is
a free $R$-module of rank $n+(n-m)(n-1)$ with basis:
$$
I,C,\dots,C^{n-1},D,CD,\dots,C^{n-m-1}D,\dots,D^{n-1},CD^{n-1},\dots,C^{n-m-1}D^{n-1}.
$$
\end{thm}

\noindent{\sl Proof.} Suppose that
\begin{equation}
\label{yu} p_0(C)I+p_1(C)D+\cdots+p_{n-1}(C)D^{n-1}=0,
\end{equation}
with $p_i\in R[X]$. We need to show that $f$ divides $p_0$ and
that $h$ divides $p_1,...,p_{n-1}$. Clearly
$$
p_1(C)D+\cdots+p_{n-1}(C)D^{n-1}=-p_0(C)I,
$$
so $p_1(C)D+\cdots+p_{n-1}(C)D^{n-1}$ commutes with $C$, which
means
\begin{equation}
\label{crucial}
p_1(C)(CD-DC)+p_2(C)(CD^2-D^2C)+\cdots+p_{n-1}(C)(CD^{n-1}-D^{n-1}C)=0.
\end{equation}
Now
$$
(CD-DC)e_1=\cdots=(CD^{n-2}-D^{n-2}C)e_1=0,
$$
so
$$
0=p_{n-1}(C)(CD^{n-1}-D^{n-1}C)e_1=p_{n-1}(C)[s].
$$
By Lemma \ref{columnspq}, $f$ divides $p_{n-1}s$ and hence
$p_{n-1}g$. It follows that $h=f/\gcd(f,g)$ divides $p_{n-1}$.
Thus by Lemma \ref{h} the last summand of (\ref{crucial}) is 0 and
can be eliminated. Proceeding like this with $e_2,\dots, e_{n-1}$
we see that $h$ divides $p_{n-2},p_{n-3},\dots,p_1$ and all these
terms can be eliminated from (\ref{crucial}). Going back to
(\ref{yu}) shows that $f$ must divide $p_0$. $\quad \blacksquare$

\begin{thm}
\label{preF} Let $R$ be a unique factorization domain and set
$h=f/\gcd(f,g)$. Let the polynomials $P_1,\dots,P_{n-1}\in R[X,Y]$
be defined as in section \ref{pj} or, more generally, be arbitrary
while satisfying (\ref{polyecu}). Then the algebra $R\langle
C,D\rangle$ has presentation
$$
\langle X,Y\,|\, f(X)=0,\, g(Y)=0,\, h(X)(X-Y)=0,\,
Y^jX=P_j(X,Y),\quad j=1,\dots,n\rangle.
$$
\end{thm}

\noindent{\sl Proof.} The proof of Theorem \ref{thA}  works as
well, except that the relation $h(A)(A-B)=0$ allows $R[A,B]$ to be
spanned by the reduced list of $n+(n-m)(n-1)$ matrices:
$$
I,A,\dots,A^{n-1},
B,AB,\dots,A^{n-m-1}B,\dots,B^{n-1},AB^{n-1},\dots,A^{n-m-1}B^{n-1}.
$$
As their images under $\Omega$ are linearly independent by theorem
\ref{basisfg}, the result follows.$\quad \blacksquare$

\section{A determinantal identity}
\label{D}

The following remarkable identity is valid for any commutative
ring $R$ with identity.

\begin{thm}
\label{thmD}
 Let the columns of $M_{f, g}\in
M_{n^2}(R)$ be the coordinates of $C^iD^j$, with $0\leq i,j\leq
n-1$, relative to the canonical basis of $M_n(R)$ formed by all basic matrices
$E^{kl}$, where $1\leq k,l\leq n$, and the lists of matrices
$C^iD^j$ and $E^{kl}$ are ordered as indicated below. Let
$M(f,g)=\mathrm{det}\, M_{f, g}$. Then $M(f,g)=R(f,g)^{n-1}$.
\end{thm}

\noindent{\sl Proof.} We order the matrices $C^iD^j$ in the
following manner:
$$D^{n-1}, CD^{n-1}, \dots , C^{n-1}D^{n-1}, \,\,D^{n-2}, CD^{n-2},
\dots , C^{n-1}D^{n-2}, \,\,\dots , \,\,I, C , \dots , C^{n-1}.$$
The basic matrices $E^{kl}$ are ordered first by column and then
by row as follows:
$$
E^{11},E^{21},\dots,E^{n1},\dots, E^{1n},E^{2n},\dots, E^{nn}.
$$
The proof consists of a sequence of reductive steps.

\begin{enumerate}

\item[(1)] Let $a\mapsto a'$ be a ring homomorphism $R\to R'$. Let
$p\to p'$ and $A\to A'$ stand for corresponding ring homomorphisms
$R[X]\to R'[X]$ and $M_n(R)\to M_n(R')$. Then
$M(f,g)=R(f,g)^{n-1}$ implies $M(f',g')=R(f',g')^{n-1}$.

This follows from the fact that $M(f,g)$ and $R(f,g)$ are defined
in such a way as to be compatible with the above ring
homomorphisms.

\item[(2)] If $R$ is an integral domain then $M(f,g)=0$ if and
only if $R(f,g)=0$.

Indeed, $M(f,g)=0$ means that the matrices $C^iD^j$ are linearly
dependent over the field of fractions of $R$, which is equivalent
to $R(f,g)=0$ by Theorem \ref{basisfg}.

\item[(3) ]  $M(f,g)$ belongs to a prime ideal $P$ of $R$ if and
only if $R(f,g)$ belongs to $P$.

This follows from (2) by using (1) with the ring homomorphism
$R\to R/P$.

\item[(4) ]  If $R$ is a unique factorization domain then $M(f,g)$
and $R(f,g)$ are both zero, both a unit, or both share
the same irreducible factors in their prime factorization.

This follows from (3).

\item[(5) ]  It suffices to prove the result for the ring
$S=\mathbb Z[Y_1,\dots,Y_{n},Z_1,\dots,Z_{n}]$.

Given $f'=a_0+\cdots+a_{n-1}X^{n-1}+X^n$ and
$g'=b_0+\cdots+b_{n-1}X^{n-1}+X^n$ in $R[X]$ we consider the ring
homomorphism $S\to R$ that restricts to the canonical map $\mathbb
Z\to R$, and sends $Y_1,\dots,Y_{n}$ to $a_0, \dots,a_{n-1}$ and
$Z_1,\dots,Z_{n}$ to $b_0, \dots,b_{n-1}$. Now use (1).

\item[(6) ] It suffices to prove the result for the field $\C$ of
complex numbers.

Clearly, to prove the result for an integral domain it is
sufficient to prove it for any field extension of its field of
fractions. In our case, $\mathbb C$ is an extension of the field
of fractions of $\Z[Y_1,\dots,Y_{n},Z_1,\dots,Z_{n}]$, so our
claim follows from (5).

\item[(7) ]  It suffices to prove the result for the ring
$S=\mathbb Z[Y_1,\dots,Y_{n},Z_1,\dots,Z_{n}]$ and the polynomials
$f=(X-Y_1)\cdots (X-Y_n)$ and $g=(X-Z_1)\cdots (X-Z_n)$.

Let $f',g'\in \mathbb C[X]$ be monic of degree $n$. Then
$f'=(X-a_1)\cdots(X-a_n)$ and $g'=(X-b_1)\cdots(X-b_n)$ for some
complex numbers $a_i,b_j$. First use (1) to derive the result for
$f'$ and $g'$ from the one for $f$ and $g$. Then apply (6).

\end{enumerate}

\medskip

We will now show that indeed $M(f,g)=R(f,g)^{n-1}$ for
$f=(X-Y_1)\cdots (X-Y_n)$ and $g=(X-Z_1)\cdots (X-Z_n)$ in $S[X]$.
This will complete the proof.

We have $R(f,g)=\Pi(Y_i-Z_j)$, with $1\leq i,j\leq n$, which is a
product of $n^2$ non-associate prime elements in the unique
factorization domain $S$. By (4) these are the prime factors of
$M(f,g)$. In particular, $R(f,g)$ divides $M(f,g)$.

Let $\sigma$ and $\tau$ be permutations of $1,\dots,n$. Let
$\Omega$ be the automorphism of $S$ corresponding to them via
$Y_i\mapsto Y_{\sigma(i)}$ and $Z_j\mapsto Z_{\tau(j)}$. This
naturally extends to automorphisms of $S[X]$ and $M_n(S)$, also
denoted by $\Omega$. As $f$ and $g$ are $\Omega$-invariant, so are
$M_{f,g}$ and $M(f,g)$.

Now if $Y_i-Z_j$ has multiplicity $a_{ij}$ in $M(f,g)$ then
$Y_{\sigma(i)}-Z_{\tau(j)}$ will have multiplicity $a_{ij}$ in
$\Omega (M(f,g))$. Since $\Omega (M(f,g))=M(f,g)$, it follows that
all prime factors of $M(f,g)$ have the same multiplicity, say
$m\geq 1$. Since the only units in $S$ are 1 and -1, we see that
$$
M(f,g)=\epsilon R(f,g)^m,\quad \epsilon\in\{1,-1\}.$$ Let
$T=\Z[Z_1,\dots,Z_n]$ and let $p\in T[X]$ be the generic
polynomial
$$
p=(X-Z_1)\cdots(X-Z_n).
$$
Then
$$
 R(f,g)=p(Y_1)\cdots p(Y_n)\in T[Y_1,\dots,Y_n].
$$

From these equations we see that the total degree of $R(f,g)$ is
$n^2$ and  the only monomial of such a degree in $R(f,g)$ is
$Y_1^n\cdots Y_n^n$, which appears with coefficient 1. Therefore
$M(f,g)$ has degree $n^2m$ and the only monomial of that degree in
$M(f,g)$ is $(Y_1\cdots Y_n)^{nm}$, which appears with coefficient
$\epsilon$. Substituting all $Z_1,\dots,Z_n$ by 0 yields
$$
M(f,X^n)=\epsilon R(f,X^n)^m,
$$
where
$$
R(f,X^n)=(Y_1\cdots Y_n)^n=(\mathrm{det}\, C)^n.
$$

We are thus reduced to proving that $M(f, X^n)= (\text{det}
\,C)^{n(n-1)}$. This we do now. Set $g=X^n$ and refer to the order of the matrices $C^iD^j$ and $E^{k l}$
given at the beginning of the proof. Expressing each vector
$C^iD^j$ in the canonical basis of $M_n(S)$ as the column vector
$$\begin{pmatrix}
C^iD^je_1\\
\vdots \\
C^iD^je_n
\end{pmatrix}\in S^{n^2},
$$
 we get the block decomposition $M_{f,\, g}\!=\! (C_{k\,j})$,
where the columns of $C_{k,\, j}\in M_n(S)$ are
$$C_{k,\, j} = (D^{n-j}e_k  \,\,CD^{n-j}e_k \,\,\dots  \,\,C^{n-1}D^{n-j}e_k),
\quad 1\le k, j\le n.$$ Let $0\le i \le n-1$, $1\le k, j \le n$.
Then
$$C^iD^{n-j}e_k=\begin{cases} C^{n-j+k-1}e_{i+1} &\text{ if $k \le j$}\\ 0 &\text{ otherwise.}\end{cases}$$
Therefore,
$$C_{k, j} =\begin{cases} C^{n-j+k-1} &\text{ if $k\le j$}\\ 0 &\text{ otherwise.}\end{cases}$$
In other words, we have
$$M_{f,\, g} = \begin{pmatrix} C^{n-1} & C^{n-2} & \dots  & C & I\\
0& C^{n-1} & \dots & C^2  &  C\\
\vdots & \vdots & \ddots & \vdots & \vdots\\
 0 & 0 & \dots & C^{n-1} & C^{n-2} \\
 0 & 0 & \dots & 0 & C^{n-1}
\end{pmatrix}.$$
Hence,
$$M(f, g) = (\text{det} \,C)^{n(n-1)}.\quad \blacksquare$$

\section{The index of $R\langle C,D\rangle$ in $M_n(R)$}
\label{indic}

Let $R$ be a principal ideal domain where each maximal ideal has
finite index. Any non-zero ideal $Ra$ is easily seen to have
finite index, which will be denoted by $N(a)$. As an example, we
may take $R$ to be the ring of integers of an algebraic number
field $\mathbb K$ of class number one, in which case
$N(a)=|N_{\mathbb K/\mathbb Q}(a)|$. In particular, $N(a)=|a|$
when $R=\mathbb Z$.

\begin{thm} Let $R$ be a principal ideal domain where each maximal ideal
has finite index in $R$. Then $R\langle C,D\rangle$ has maximal
rank in $M_n(R)$ if and only if $R(f,g)\neq 0$, in which case
$[M_n(R):R\langle C,D\rangle]=N(R(f,g))^{n-1}$.
\end{thm}

\noindent{\sl Proof.} Let $R^*$ be the monoid of non-zero elements
of $R$ and write $\mathbb N$ for the monoid of natural numbers. By
hypothesis each maximal ideal $Rp$ has finite index, denoted by
$N(p)$.

If $a\in R$ is not zero or a unit then $a=p_1^{a_1}\cdots
p_m^{a_m}$, where the $p_i$ are non-associate primes in $R$ and
$a_i\geq 1$. Using the Chinese Remainder Theorem and the fact that
$p^iR/p^{i+1}R$ is a one-dimensional vector space over $R/p$ for
every prime $p$, it follows at once that $Ra$ also has finite
index, say $N(a)$, in $R$, where $N(a)=N(p_1)^{a_1}\cdots
N(p_m)^{a_m}$. Thus $N:R^*\to\mathbb N$ is a homomorphism of
monoids whose kernel is the unit group of $R$.

We have the free $R$-module of rank $M_n(R)$ of rank $n^2$ and its
submodule $R\langle C,D\rangle$, which is free of rank $\leq n^2$.
By Corollary \ref{iguales} the matrices $C^iD^j$, with $0\leq
i,j\leq n-1$, span $R\langle C,D\rangle$. The matrix expressing
the coordinates of these generators in the basis of $M_n(R)$
formed by all $E^{ij}$ is the matrix $M_{f, g}$ of Theorem
\ref{thmD}. Let $a_1,\dots , a_{n^2}$ be the invariant factors of
$M_{f, g}$. Then $M_n(R)$ has a basis $u_1,\dots,u_{n^2}$ such
that $a_1u_1,\dots,a_{n^2}u_{n^2}$ span $R\langle C, D\rangle$.
Hence $R\langle C, D\rangle$ has rank $n^2$ if and only if
$M(f,g)=a_1\cdots a_{n^2}\neq 0$. Since $M_n(R)/R\langle C,
D\rangle\cong R/Ra_1\times\cdots \times R/Ra_{n^2}$ as
$R$-modules, if all $a_1,\dots,a_{n^2}$ are non-zero then
$[M_n(R):R\langle C, D\rangle]=N(a_1\cdots a_{n^2})=N(M(f, g)).$
Now apply Theorem \ref{thmD}.$\quad \blacksquare$

\bigskip

\noindent{\bf {\Large Acknowledgements}}

\bigskip

The authors thank D. Stanley for useful conversations and D.
Djokovic for writing a computer program to verify that Theorem
\ref{thmD} was indeed true when $n=3,4$.

\bigskip

\noindent{\bf{\Large References}}

\medskip

\small

\begin{enumerate}
\item[[A]$\!\!\!$] G. Agnarsson, On a class of presentations of
matrix algebras. Comm. Algebra 24 (1996), 4331-4338.

\item[[AM]$\!\!\!$]  M.F. Atiyah and I.G. Macdonald, Introduction
to commutative algebra, Addison-Wesley, 1969.

\item[[AMR]$\!\!\!$] G. Agnarsson, S.A. Amitsur and J.C.
Robson, Recognition of matrix rings II, Israel J. Math. 96 (1996),
1-13.

\item[[CR]$\!\!\!$] C.W. Curtis, and I. Reiner, Representation theory of finite groups and associative algebras,
Interscience, 1962.

\item[ [C1]$\!\!\!$] A.W. Chatters, Representation of tiled matrix
rings as full matrix rings, Math. Proc. Cambridge Philos. Soc. 105
(1989), 67-72.

\item[[C2]$\!\!\!$] A.W. Chatters, Matrices, idealisers, and
integer quaternions, J. Algebra 150 (1992), 45-56.

\item[[GS]$\!\!\!$] N.H. Guersenzvaig and F. Szechtman,
A closed formula for the product in simple integral extensions,
Linear Algebra Appl. 430 (2009), 2464-2466.

\item[ [LRS]$\!\!\!$] L.S. Levy, J.C. Robson and J.T. Stafford,
Hidden matrices, Proc. London Math. Soc. (3) 69 (1994), 277-308.

\item[ [R]$\!\!\!$] J.C. Robson, Recognition of matrix rings,
Comm. Algebra 19 (1991), 2113-2124.
\end{enumerate}

\end{document}